# Treating the Future Equally: Solving Undiscounted Infinite Horizon Optimization Problems

Dapeng CAI [1] and Gyoshin NITTA [2]

[1] *Institute for Advanced Research, Nagoya University, Furo-cho, Chikusa-ku, Nagoya, 464-8601, Japan;* [2] *Department of Mathematics, Faculty of Education, Mie University, Kurimamachiya 1577, Tsu, 514-8507, Japan*

**Abstract**

Infinite horizon optimization problems accompany two perplexities. First, the infinite series of utility sequences may diverge. Second, boundary conditions at the infinite terminal time may not be rigorously expressed. In this paper, we show that under two fairly general conditions, the limit of the solution to the undiscounted finite horizon problem is optimal among feasible paths for the undiscounted infinite horizon problem, in the sense of the overtaking criterion. Applied to a simple Ramsey model, we show that the derived path contains intriguing properties. We also comprehend the legitimacy of the derived paths by addressing the perplexities with non-standard arguments.

*True end is not in the reaching of the limit, but in a completion which is limitless.*

- *The Crescent Moon*, Rabindranath Tagore

## 1. Introduction

Increasingly economists are being asked to evaluate and compare policies whose effects can be expected to spread out into the far distant future. Prominent examples are abundant: global climate change, radioactive waste disposal, loss of biodiversity, groundwater pollution, minerals depletion, and many others (Nordhaus, 1994; Weitzman, 1998; 1999)

However, such problems are "genuinely deep and difficult", as "neat and convincing general answer[s]" may not exist (Solow, 1999, p. vii). With an infinite horizon, the infinite series of utility sequences in general will diverge. Consequently, one ends up at comparing infinity with infinity, which has been impossible within the real number field. As a compromise, economists have been assuming pure time preference, although they have long been scathing about its ethical dimensions. Indeed, discounting has dominated intertemporal economics, especially the benefit-cost analysis, "more for lack of convincing alternatives than because of the conviction it inspires (Heal, 1998, p. 12)." The sense of the unease with discounting starts with Ramsey, who commented that discounting "is ethically indefensible and arises merely from the weakness of the imagination (1928, p. 543)." Pigou stated that pure time preference "implies … our telescopic faculty is defective (1932, p. 25)." Harrod also expressed his reservations by saying that it is a "human infirmity" and "a polite expression for rapacity and the conquest of reason by passion (1948, p. 40)." Solow held that "in solemn conclave assembled, so to speak, we ought to act as if the social rate of time preference were zero (1974, p. 9)." Cline (1992) and Anand and Sen (2000) also take similar stands. Today, the consensus on this issue seems to be that "we [should] treat the welfare of future generations on a par with our own" and uncertainties about future prospects, except for the case of "uncertainty about existence", "are not reasons for discounting (Stern, 2007, p. 35-36)." Discounting raises moral and logical difficulties because a positive discount rate generates a fundamental asymmetry between the treatments of the present and future generations, which is most troublesome when applied to environmental matters such as climate change that happen in the deep future (Heal, 1998). Because of the very counterintuitive "relentless force of compound interest (Weitzman, 1998, p. 202)", "…discounting at anything like market interest rates implies conclusions that common sense can not accept (Solow, 1999, p. viii)." Indeed, "even very large damages, if they will happen 200 years from now, discount back to peanuts (Solow, *ibid.*)." Hence, "it would not make much difference in the relevant unit of today's standard of living whether climate change from global warming comes 200 years from now or 400 years from now (Weitzman, *ibid.*)."

As pointed in Koopmans (1960), the technical difficulties of not discounting lies in the fact that "there is not enough room in the set of real numbers to accommodate and label numerically all the different satisfaction levels that may occur in relation to consumption programs for an infinite future (p. 288)", when there is "a preference for postponing satisfaction, or even neutrality toward timing (*ibid.*)." Obviously, to overcome these difficulties, one needs to extend the real number field $\mathbb{R}$ to one that could accommodate and distinguish the many "levels" of infinities. Such a number field is provided by the theory of non-standard analysis (Robinson, 1966). Non-standard analysis has been applied to probability theory and mathematical economics (Anderson, 1991).[1] Rubio has applied non-standard analysis to optimization theory (Rubio, 1994, 2000).[2] However, by far such application has been incomplete and not readily applicable to most problems in economic dynamics.

---

[1] Non-standard methods have been used to examine large economies and the probability theory. An excellent survey on the economic applications of the non-standard analysis is available in Anderson (1991), which provides "a careful development of non-standard methodology in sufficient detail to allow the reader to use it in diverse areas in mathematical economics" (p. 2147).

[2] Rubio (2000, Appendix) also provides a brief, yet extremely helpful introduction to the non-standard analysis.

In the present paper, we consider the undiscounted version of an otherwise standard one-sector Ramsey model of economic growth. The derivation of our main results does not necessarily require the knowledge of the non-standard analysis. For our purpose, the use of the overtaking criterion would be enough, as it does enable a comparison of infinities in the real number field.[3] We show that under two fairly general conditions in the field of real numbers, the limit of the solution to the undiscounted finite horizon problem is optimal among all feasible paths for the undiscounted infinite horizon problem, in the sense of the overtaking criterion. The conjecture that the limit of the solutions to the finite horizon problem is the unique solution to the infinite horizon problem has been around for a while.[4] By demonstrating that the conjecture is correct, we have presented a fruitful way of analyzing undiscounted infinite horizon optimization problems. One would be able to explicitly analyze the problem, in much the same way as one would treat the undiscounted finite horizon problem, as long as the model in consideration satisfies the two conditions.

Nevertheless, a slightly touch of the non-standard analysis would enable one to better comprehend the legitimacy of the conjecture, which mandates a correct presentation of the problem. This has been impossible in the real number field: as "∞+1" is indefinable, the boundary condition that "goods will not be wasted" at the infinite terminal time cannot be appropriately presented. By resorting to the concepts of the non-standard analysis, we are able to reformulate the optimization question in a mathematically correct form and show that "the limit of the solutions to the finite horizon problem" is the unique "projection" of the unique hyper-real optimal solution on the real number field. Hence we have shown that the conjecture is meaningful, as "the limit of the solutions to the finite horizon problem" is the standard part of the unique optimal path in the hyper-real number field that satisfies the boundary condition at the infinite terminal time. On the other hand, the application of our approach does not require the knowledge of the non-standard analysis. For most examples, it may be possible to use a numerical approach to check the two conditions and to compute explicit solutions.

---

[3] However, different from the non-standard analysis, such a comparison is indirect, and as a result, the legitimacy of the derived path cannot be established with the overtaking criterion.

[4] The case with a discount factor has been examined by using the recursive method in for example, Stokey and Lucas (1989).

In the following section we begin by considering a simple undiscounted one-sector finite horizon Ramsey model of economic growth. By resorting to the arguments in the non-standard analysis, we extend real numbers to hyper-real numbers and present a mathematically correct infinite horizon extension of the finite horizon problem. We then demonstrate that there exists a unique optimal solution to this infinite horizon problem in the hyper-real number field, with the "projection" of which on the real number field being the limit of the solution to the undiscounted finite horizon problem. We proceed to prove that under two fairly general real conditions, such a path is optimal among all feasible paths in the real number field, in the sense of the overtaking criterion. We then apply our approach to a particular parametric example of the Ramsey model. We show that as ethical ramifications of discounting large future costs and benefits into small present values are nonexistent, the optimal paths exhibit some intuitively appealing properties.

## 2. The Model

### 2.1 Two Perplexities of the Infinite Horizon Optimization Problem

We consider an economy that is composed of many identical households, each forming an immortal extended family. In each period, a representative household invests $k(t)$ at the beginning of the period $t$ to produce $f(k(t))$ amount of output. The production process is postulated as follows: $f: \mathbb{R}_+ \to \mathbb{R}_+$, $f$ is continuous, strictly increasing, weakly concave, and continuous differentiable, with $f(0) = 0$, and for some $\bar{x} > 0: x \leq f(x) \leq \bar{x}$, if $0 \leq x \leq \bar{x}$; and $f(x) < x$, if $x > \bar{x}$. In each period, output is divided between current consumption $c(t)$, and investment, $k(t+1)$, with $c(t) \geq 0$, and $k(t+1) \geq 0$. Given the planning horizon $T \in [0, \infty)$, the criterion for a social planner to judge the welfare of the representative household takes the form

$$\max_{c(t)} \sum_{t=0}^{T} U(c(t)), \tag{1}$$

where the instantaneous utility function $U: \mathbb{R}_+ \to \mathbb{R}$ is continuous, strictly increasing, strictly concave, and continuous differentiable. The household chooses the path of $c(t)$ that maximizes (1), which is subject to the budget constraint: $k(t+1) = f(k(t)) - c(t)$,

the initial capital stock: $k(0) \equiv k_0 > 0$, and

the boundary condition: $k(T+1) = 0$. (2)

A unique solution to problem (1), subject to (2), $\{c_T(t), k_T(t)\}_{t=0}^{T}$, can be found readily by applying the method of Lagrange, with the first-order conditions given as[5]

$$U'(c_T(t)) = \lambda_T(t), \tag{3}$$

$$\lambda_T(t) f'(k_T(t)) = \lambda_T(t-1). \tag{4}$$

We then have

$$f'(k_T(t)) U'(f(k_T(t)) - k_T(t+1)) = U'(f(k_T(t-1)) - k_T(t)). \tag{5}$$

---

[5] The subscripts denote the length of the planning horizon.

We proceed to extend the planning horizon of this problem to infinity. Two perplexities arise immediately. First of all, the infinite series $\sum_{t=0}^{\infty} U(c(t))$ will in general diverge and the maximization of which may be meaningless in the real number world. Second, it would also be imperative to present an appropriate analogue of the boundary condition (2) for the infinite-horizon case. The naïve extension of which would be $k(\infty+1) = 0$. However, as $\infty$ has not been defined in $\mathbb{R}$, $k(\infty+1) = 0$ is indefinable, and we can only state that "goods will not be wasted" in the context of real number field, as in Stokey and Lucas (1989). It seems that the boundary condition can be stated instead as $\lim_{T \to \infty} k(T+1) = 0$, however, which also implies $\lim_{T \to \infty} k(T) = 0$. Since the capital stock $k$ at the infinite time is not specified as zero, we have arrived at a contradiction and it is inappropriate to state the boundary condition with limiting concepts in $\mathbb{R}$.

## 2.2  A Reformulation of the Problem

Instead, we extend the field of real numbers so that different "levels of infinity" can be accommodated and we can distinguish "$\infty$" and "$\infty+1$". Following Anderson (1999, pp. 2150-2151), we extend real numbers into hyper-real numbers. We define hyper-real number as $^*\mathbb{R} \equiv \mathbb{R}^{\mathbb{N}}/u$, where $\mathbb{N}$ is a set of natural numbers, $u$ is a free ultrafilter that is a maximal filter containing the *Fréchet filter*. The elements of $^*\mathbb{R}$ are represented as sequences and are denoted as $[<a_n>]$, where $a_n \in \mathbb{R}$. Instead of $T$, we consider an infinite star finite number $\tilde{T} \equiv [<T_n>]$. We extend $c(t)$ to $C(t) \equiv [<C_n(t)>]$, and $k(t)$ to $K(t) \equiv [<K_n(t)>]$. Moreover, $U : \mathbb{R}_+ \to \mathbb{R}$ and $f : \mathbb{R}_+ \to \mathbb{R}_+$ are extended to $^*U : {^*\mathbb{R}_+} \to {^*\mathbb{R}}$ and $^*f : {^*\mathbb{R}_+} \to {^*\mathbb{R}_+}$, respectively. Accordingly, $\sum_{t=0}^{T} U(c(t))$ is extended to $^*\sum_{[<t_n>]=0}^{\tilde{T}} {^*U(C(t))}$, with its elements represented by $\left[ < \sum_{t=0}^{T_n} U(C_n(t)) > \right]$. We also denote

$$[<a_n>] \geq [<b_n>] \text{ if } \{n \in \mathbb{N} : a_n \geq b_n\} \in u,$$

for $[<a_n>], [<b_n>] \in {^*\mathbb{R}}$.

Hence, the maximization problem is reformulated as

$$\max {^*\sum_{[<t_n>]=0}^{\tilde{T}}} {^*U(C(t))}, \qquad *(1)$$

subject to

$$C([<t_n>]) = {^*f(K([<t_n>]))} - K([<t_n>]+1),$$

$$0 \leq K([<t_n>]+1) \leq {^*f(K([<t_n>]))}, \quad t_n = 0,1,\ldots,T_n,$$

given $K(0) > 0$, and,

the fact that goods will not be wasted: $K(\tilde{T}+1)=0$. *(2)

Next, we show how to solve the optimization problem in the hyper-real number field.

Let $[<\lambda_n>]:\{1,2,\ldots,\tilde{T}\} \to {}^*\mathbb{R}$, the first-order conditions are familiar:

$${}^*U'(C([<t_n>])) = [<\lambda_n>]([<t_n>]),$$ *(3)

$$[<\lambda_n>]([<t_n>])\,{}^*f'(K([<t_n>])) = [<\lambda_n>]([<t_n>]-1).$$ *(4)

Combining *(3) and *(4), we have

$${}^*f'(K[<t_n>])\,{}^*U'({}^*f(K([<t_n>])) - K([<t_n>]+1)) = {}^*U'({}^*f(K([<t_n>]-1)) - K([<t_n>])),$$

where $t = 1, 2, \ldots, \tilde{T}$. *(5)

Equation *(5) is a second-order difference equation in $K(t)$; hence it has a two-parameter family of solutions. The unique optimum of interest is the one solution in this family that in addition satisfies the two boundary conditions

$K(0) > 0$, $K(\tilde{T}+1) = 0$ given.

For each index $n$, equation *(5) is equivalent to

$$f'(K_n(t_n))U'[f(K_n(t_n)) - K_n(t_n+1)] = U'[f(K_n(t_n-1)) - K_n(t_n)],$$

where $t_n = 1, 2, \ldots, T_n$, *(5')

and an unique solution $\left(C_{\tilde{T}}(t), K_{\tilde{T}}(t)\right)$ can be obtained by applying the following boundary conditions:

$$K_n(0) > 0, \quad K_n(T_n + 1) = 0.$$

The variables in equations *(3), *(4), and *(5) all depend on $\tilde{T}$. In particular, $[<\lambda_n>]([<t_n>])$ stands for the *shadow price for capital at time* $[<t_n>]$, *given the planning horizon* $\tilde{T}$. In other words, the value of $\lambda$ at each time changes with the length of the planning horizon. Following the arguments of the non-standard analysis, we then have

**Theorem 1.** There exists a unique optimal solution $\{C(t), K(t)\}_{t=0}^{\tilde{T}}$ to problem *(1), subject to *(2).

An intuitive interpretation of the theorem is as follows: under our extension of $\mathbb{R}$ to $^*\mathbb{R}$, the relation between $T$ and $T+1$ has been preserved to $\tilde{T}$ and $\tilde{T}+1$. In other words, the set of sequences $\{K(t)\}_{t=0}^{\tilde{T}}$ satisfying *(2) is a *closed, *bounded, and *convex subset of $^*\mathbb{R}$, and the objective function *(1) is continuous and strictly concave. Therefore, there exists a unique solution $\{C(t), K(t)\}_{t=0}^{\tilde{T}}$ to *(1) in the hyper-real number field.

2.3 The Projection of the Hyper-real Optimal Solutions on the Real Number Field

Denote $\widetilde{st}$ to be the standard mapping from $^*\mathbb{R} \to \mathbb{R} \cup \{\pm\infty\}$. $\widetilde{st}$ projects the unique optimal path in the hyper-real number field to the real number field. Two conversions accompany this projection: $\widetilde{T}$ is converted to $\infty$, and the infinitesimals are neglected. As is widely known in the non-standard analysis, if there exist $\lim_{T \to \infty} c_T(t)$ and $\lim_{T \to \infty} k_T(t)$, we then have

$$\widetilde{st}\left(C_{\widetilde{T}}(t)\right) = \lim_{T \to \infty} c_T(t) \equiv c^\circ(t), \widetilde{st}\left(K_{\widetilde{T}}(t)\right) = \lim_{T \to \infty} k_T(t) \equiv k^\circ(t),$$

where $t \in \mathbb{R}$.

Hence, from the non-standard optimal conditions, we have specified a unique path $\left(\widetilde{st}\left(C_{\widetilde{T}}(t)\right), \widetilde{st}\left(K_{\widetilde{T}}(t)\right)\right)$ in $\mathbb{R}$, although we still need to verify whether such a path is indeed optimal in $\mathbb{R}$.

### 2.4 The Overtaking Criterion

In what follows, we show that under two fairly general conditions, the optimum obtained by applying Theorem 1 is indeed optimal among all feasible paths in $\mathbb{R}$. For feasible paths that satisfy the following conditions:

the budget constraint: $k(t+1) = f(k(t)) - c(t)$,

the initial capital stock: $k(0) \equiv k_0 > 0$, and the boundary condition that no goods should be wasted,

we convert them by changing the values of both consumption and capital at a newly chosen "terminal date", $T$, so that both are 0 at $T+1$:[6]

$$\text{for all } T, \quad \tilde{c}(t) = \begin{cases} c(t), 0 < t < T-1, \\ f(k(T)), t = T, \\ 0, t \geq T+1, \end{cases} \quad \tilde{k}(t) = \begin{cases} k(t), 0 < t \leq T, \\ 0, t \geq T+1. \end{cases}$$

This endeavour enables us to explicitly incorporate the boundary condition that no goods should be wasted. For two such converted feasible paths $(k_1, c_1)$ and $(k_2, c_2)$, we define the order relation "$\leq$" (i.e., $(c_2, k_2)$ catches up to (overtakes) $(c_1, k_1)$) as follows:

**Definition.** For $U(c_1(t)) < 0$, $U(c_2(t)) < 0$, all $t$,

---

[6] Our formulation of the overtaking criterion

$$a(t) > b(t) \quad \text{if} \quad \lim_{T \to \infty} \left( \sum_{t=0}^{T-1} (a(t) - b(t)) + \tilde{a}(T) - \tilde{b}(T) \right) > 0,$$

appears to be different from those commonly used in the literature,

$$a(t) > b(t) \quad \text{if} \quad \lim_{T \to \infty} \sum_{t=0}^{T} (a(t) - b(t)) > 0.$$

However, the two are identical when $\lim_{T \to \infty} (\tilde{a}(T) - \tilde{b}(T))$ and $\lim_{T \to \infty} (a(T) - b(T))$ are finite, whereas $\lim_{T \to \infty} \sum_{t=0}^{T} (a(t) - b(t))$ is infinite.

$$(k_1, c_1) \leq (k_2, c_2) \ if \ \lim_{T \to \infty} \frac{\sum_{t=0}^{T-1} U(c_1(t)) + U(f(k_1(T)))}{\sum_{t=0}^{T-1} U(c_2(t)) + U(f(k_2(T)))} \geq 1.^{7,\ 8}$$

Following Brock's (1970) notion of weak maximality, the optimality criterion in this paper is defined as follows: a feasible path $(k, c)$ is optimal if it catches up to (overtakes) any other feasible path $(k_1, c_1)$.[9] In Theorem 2, we consider the case for $U(c(t)) < 0$ and $U(c_1(t)) < 0$, all $t$. In that case, the optimality criterion is defined as

---

[7] Our overtaking criterion also differs from those in the literature in that it is in the form of

$\sum_{t=0}^{T} a(t) / \sum_{t=0}^{T} b(t)$, rather than $\sum_{t=0}^{T} a(t) - \sum_{t=0}^{T} b(t)$. These two forms are equivalent as

$\sum_{t=0}^{T} a(t) - \sum_{t=0}^{T} b(t) = \sum_{t=0}^{T} b(t) \left( \sum_{t=0}^{T} a(t) / \sum_{t=0}^{T} b(t) - 1 \right)$, if there exists $\lim_{T \to \infty} \sum_{t=0}^{T} b(t) > 0$. We use

the former form for the clarity of presentation concerning the proof of condition (ii) in Theorem 2. In other words, our criterion is $\infty / \infty$ type, whereas the one in the literature is $(\infty - \infty)$ type.

[8] $\tilde{T} = \left[ \tilde{T}_n \right]_{n=1}^{\infty}$ denotes the time at which the system is evaluated. As $n$ can be any positive

increasing sequence, the time span of evaluation may differ. For $\tilde{T}$, there exists an overtaking criterion $\lim_{n \to \infty} \sum_{t=1}^{\tilde{T}_n} (a(t) - b(t)) > 0$. Specially, when $\tilde{T} = [1, 2, 3, \cdots]$, the criterion is the usual standard overtaking criterion.

[9] According to the non-standard argument, the standard parts of the two derived paths are exactly their original paths when $0 < t < \infty$, respectively. Moreover, our optimal programs are not restricted to optimal stationary programs.

$$(k_1, c_1) \leq (k, c) \text{ if } \lim_{T \to \infty} \frac{\sum_{t=0}^{T-1} U(c_1(t)) + U(f(k_1(T)))}{\sum_{t=0}^{T-1} U(c(t)) + U(f(k(T)))} \geq 1.$$

The results for the case for $U(c(t)) > 0$ and $U(c_1(t)) > 0$ can be obtained similarly.

## 2.5 The Proof of the Conjecture

To establish the legitimacy of the conjecture, we would first need the following lemma:

**Lemma 1.** Let $a_T$ and $b_T$, $T \in [0, \infty)$, be two sequences. If $\lim_{T \to \infty} a_T > 0$, $\lim_{T \to \infty} b_T > 0$, then $\lim_{T \to \infty} (ab)_T = \lim_{T \to \infty} a_T \cdot \lim_{T \to \infty} b_T$.

**Proof.** See Appendix.

**Theorem 2.** If the following conditions (i) and (ii) are simultaneously satisfied, then for an arbitrary $(k_1, c_1)$, $(k_1, c_1) \leq (k, c)$:

Condition (i): $\lim_{T \to \infty} \sum_{t=0}^{T} U(c_T(t))$ is infinite.

Condition (ii): $\lim_{T \to \infty} \frac{\sum_{t=0}^{T-1} \left( U(c_T(t)) - U(c^\circ(t)) \right) + U(f(k_T(T))) - U(f(k^\circ(T)))}{\sum_{t=0}^{T} U(c_T(t))} = 0.$

**Proof.** For an arbitrary number $T$, as $c_T(t)$ is the optimal solution, we have

$$\frac{\sum_{t=0}^{T-1}U(c_1(t))+U(f(k_1(T)))}{\sum_{t=0}^{T}U(c_T(t))}\geq 1.$$

Also, as

$$\frac{\sum_{t=0}^{T-1}U(c_1(t))+U(f(k_1(T)))}{\sum_{t=0}^{T-1}U(c^\circ(t))+U(f(k^\circ(T)))}=$$

$$\left(\frac{\sum_{t=0}^{T-1}U(c_1(t))+U(f(k_1(T)))}{\sum_{t=0}^{T}U(c_T(t))}\right)\cdot\left(\frac{\sum_{t=0}^{T}U(c_T(t))}{\sum_{t=0}^{T-1}U(c^\circ(t))+U(f(k^\circ(T)))}\right),$$

and since $c_T(T)=f(k_T(t))$,

$$\frac{\sum_{t=0}^{T}U(c_T(t))}{\sum_{t=0}^{T-1}U(c^\circ(t))+U(f(k^\circ(T)))}$$

$$=\frac{\sum_{t=0}^{T}U(c_T(t))}{\sum_{t=0}^{T}U(c_T(t))-\left(\sum_{t=0}^{T-1}U(c_T(t))-\sum_{t=0}^{T-1}U(c^\circ(t))\right)-U(c_T(T))+U(f(k^\circ(T)))}$$

$$=\frac{1}{1-\dfrac{\sum_{t=0}^{T-1}\bigl(U(c_T(t))-U(c^\circ(t))\bigr)+U(f(k_T(T)))-U(f(k^\circ(T)))}{\sum_{t=0}^{T}U(c_T(t))}}.$$

From condition (ii) in Theorem 2, we see that $\displaystyle\lim_{T\to\infty}\frac{\sum_{t=0}^{T}U(c_T(t))}{\sum_{t=0}^{T-1}U(c^\circ(t))+U(f(k^\circ(T)))}=1$.

Hence, from Lemma 1, we have

$$\lim_{T\to\infty} \frac{\sum_{t=0}^{T-1} U(c_1(t)) + U(f(k_1(T)))}{\sum_{t=0}^{T-1} U(c^\circ(t)) + U(f(k^\circ(T)))}$$

$$= \lim_{T\to\infty} \left( \frac{\sum_{t=0}^{T-1} U(c_1(t)) + U(f(k_1(T)))}{\sum_{t=0}^{T} U(c_T(t))} \right) \cdot \lim_{T\to\infty} \left( \frac{\sum_{t=0}^{T} U(c_T(t))}{\sum_{t=0}^{T-1} U(c^\circ(t)) + U(f(k^\circ(T)))} \right) \geq 1.$$

<div align="right">*Q.E.D.*</div>

Hence we have shown that the path obtained from Theorem 1 is *indeed* optimal among all feasible paths in $\mathbb{R}$, if the two conditions listed in Theorem 2 are simultaneously satisfied.

If Condition (i) is not satisfied, *i.e.*, if the maximand is finite, we can always resort to the standard Lagrange method to find the optimum. On the other hand, Condition (ii) requires that the loss accompanying the limiting operations is negligible as compared to the value of the maximand.

In the following section, we study a particular parametric example of the simple Ramsey model. For most other parametric examples, it may not be possible to explicitly check the two conditions in Theorem 2 and study the resultant paths. In such cases a numerical approach can be used to check the two conditions and to compute explicit solutions.

### 3. An Example

In this section, we first show that for a particular parametric example of the simple Ramsey model, the derived path is indeed the unique optimum among all feasible paths in $\mathbb{R}$. We also show that the path exhibit some intriguing properties.

### 3.1 A Parametric Example

We consider the following social planner's problem, in which the social planner treats the future equally by assuming a zero rate of time preference. The planner's objective is

$$\max_{c(t)} \sum_{t=0}^{T} \ln(c(t)), \text{ where } T \in [0, \infty), \tag{6}$$

subject to

$$c(t) + k(t+1) = f(k(t)),$$

$$\text{where } f(k(t)) = k(t)^\alpha, \ 0 < \alpha < 1, \ k(0) \text{ given,}$$

$$0 < k(t) < 1, \ k(T+1) = 0. \tag{7}$$

We extend $T \to \tilde{T}$, (6) and (7) are reformulated as

$$\max_{C(t)} {}^*\!\!\!\sum_{[<t_n>]=0}^{\tilde{T}} {}^*\ln(C(t)), \tag{*6}$$

subject to

$$C([<t_n>]) + K([<t_n>]+1) = (K([<t_n>]))^\alpha,$$

where $0 < \alpha < 1$, $K(0)$ given,

$$0 < K([<t_n>]) < 1, \quad K(\tilde{T}+1) = 0. \qquad *(7)$$

As $0 < k(t) < 1$, we have $0 < c(t) < 1$. The paths of $c(t)$ and $k(t)$ are uniquely determined given $k(0) > 0$ and $k(T+1) = 0$, with $k(t+1) = \alpha \dfrac{1-\alpha^{T-t}}{1-\alpha^{T-t+1}} k^\alpha(t)$, $t = 0, 1, \ldots, T$. Also, the optimal solution $(k_T(t), c_T(t))$ is

$$k_T(t) = \frac{\left(1-\alpha^{T-t+1}\right)\alpha^{\frac{1-\alpha^t}{1-\alpha}} k(0)^{\alpha^t}}{\left(1-\alpha^{T-t+2}\right)^{1-\alpha}\left(1-\alpha^{T-t+3}\right)^{\alpha(1-\alpha)} \cdots \left(1-\alpha^{T+1}\right)^{\alpha^{t-1}}},$$

$$c_T(t) = \frac{(1-\alpha)\alpha^{\frac{\alpha(1-\alpha^t)}{1-\alpha}} k(0)^{\alpha^t}}{\left(1-\alpha^{T-t+1}\right)^{1-\alpha}\left(1-\alpha^{T-t+2}\right)^{\alpha(1-\alpha)} \cdots \left(1-\alpha^{T}\right)^{\alpha^{t-1}(1-\alpha)}\left(1-\alpha^{T+1}\right)^{\alpha^t}}.$$

Next we extend $T \to \tilde{T}$, where $\tilde{T}$ is an infinite star finite number. Following Theorem 1, we then have

$$\widetilde{st}\left(K_{\tilde{T}}(t+1)\right) = \lim_{T \to \infty} k_T(t+1) = \alpha \left(\widetilde{st}\left(K_{\tilde{T}}(t)\right)\right)^\alpha, \quad t \in [0, \infty).$$

Hence,

$$\widetilde{st}\left(K_{\tilde{T}}(t)\right) = \lim_{T \to \infty} k_T(t) = (\alpha)^{1/(1-\alpha)} \left(\frac{k(0)}{(\alpha)^{1/(1-\alpha)}}\right)^{\alpha^t},$$

$$\widetilde{st}\,(C_{\tilde{T}}(t)) = \lim_{T \to \infty} c_T(t) = (1-\alpha)(\alpha)^{\alpha/(1-\alpha)} \left( \frac{k(0)}{(\alpha)^{1/(1-\alpha)}} \right)^{\alpha^{t+1}},$$

$$\widetilde{st}\,(\lambda_{\tilde{T}}(t)) = \lim_{T \to \infty} \lambda_T(t) = (1-\alpha)^{-1}(\alpha)^{-\alpha/(1-\alpha)} \left( \frac{k(0)}{(\alpha)^{1/(1-\alpha)}} \right)^{-\alpha^{t+1}}.$$

### 3.2   The Properties of the Derived Optimal Path

Denoting the derived optimal paths as $\widetilde{st}(K_{\tilde{T}}(t)) \equiv k°(t),\ \widetilde{st}(C_{\tilde{T}}(t)) \equiv c°(t),$ $\widetilde{st}(\lambda_{\tilde{T}}(t)) \equiv \lambda°(t),$ we have:

**Result 1**. The optimal saving ratio $r$ is constant over time: $r \equiv \dfrac{k°(t+1)}{f(k°(t))} = \alpha$.

Ramsey considers the continuous-time version of the problem. Under his formulation, it is shown that the optimal saving ratio, which does not depend on the productivity, is constant over time and can be very high (Ramsey, 1928).[10] Here we show that the optimal saving ratio depends on the productivity $\alpha$ (as the elasticity of output with respect to capital $\alpha$ is about the amount of output obtainable from a given unit of capital), which may lead to different estimation of the optimal saving ratio.

---

[10] Arrow concludes that "the strong ethical requirement that all generations be treated alike, itself reasonable, contradicts a very strong intuition that it is not morally acceptable to demand excessively high savings rates of any one generation, or even of every generation (Arrow, 1999, p. 16)." However, as argued in Stern, Arrow's argument is not convincing as, first, Arrow's model has not considered all aspects influencing optimum saving, and second, "his way of 'solving' the 'over-saving' complication is very ad hoc (Stern, 2007, p. 54)."

We also consider the dynamics of the model, the merits of technology ($\alpha$), the costs of not applying the optimization practice immediately, and the right values of capital at each time. Summarized as proposition 1-4, we show that the new optimal path demonstrate some intriguing properties. Moreover, in Appendix, we show that the two conditions in Theorem 2 are satisfied and the derived path is indeed optimal among all feasible paths in the real number field.

**Proposition 1.** (Convergence) $c°(t)$, $k°(t)$, and $\lambda°(t)$ converge to finite positive values as $t \to \infty$: $c_\infty \equiv \lim_{t \to \infty} c°(t) = (1-\alpha)(\alpha)^{\frac{\alpha}{1-\alpha}} > 0$, $k_\infty \equiv \lim_{t \to \infty} k°(t) = (\alpha)^{\frac{1}{1-\alpha}} > 0$, $\lambda_\infty \equiv \lim_{t \to \infty} \lambda°(t) = (1-\alpha)^{-1}(\alpha)^{\frac{-\alpha}{1-\alpha}} > 0$. When $k(0) > (\alpha)^{\frac{1}{1-\alpha}}$ $\left(k(0) < (\alpha)^{\frac{1}{1-\alpha}}\right)$, $k°(t)$ and $c°(t)$ are monotonically decreasing (increasing) in $t$, whereas $\lambda°(t)$ is monotonically increasing (decreasing) in $t$. $k°(t)$, $c°(t)$, and $\lambda°(t)$ are in steady states when $k(0) = (\alpha)^{\frac{1}{1-\alpha}}$.

Proposition 1 indicates that the economy converges to a steady state that is different from the one specified in the literature, precisely because that different from the literature, the relevant boundary conditions at $T = \infty$ have been incorporated, as an extension of the finite-horizon problem.

**Proposition 2.** (The elasticity effect) The effect of the elasticity of output with respect to capital, $\alpha$, dominates the effect of initial capital value $k(0)$ as $t \to \infty$. In the steady state, both $k_\infty$ and $c_\infty$ are decreasing in the $\alpha$.

Proposition 2 follows immediately from Proposition 1. It shows that as the system approaches the steady state, the effects of the initial capital value gradually vanishes away, whereas the effect of elasticity of output enhanced along the way, with the steady state exclusively depends on the latter. We see the higher the level of the elasticity, the higher are the stationary values in the steady state.

**Proposition 3.** (Path-dependence) It takes infinite time for $k°(t)$ and $c°(t)$ to approach the steady states. Nevertheless, the time needed depends on the initial value of $k(0)$. Formally, the difference of time needed for two initial values $k^a(0) \equiv \frac{1}{a}(\alpha)^{\frac{1}{1-\alpha}}$ and $k^b(0) \equiv \frac{1}{b}(\alpha)^{\frac{1}{1-\alpha}}$ is $\frac{\ln\left(\frac{\ln a}{\ln b}\right)}{-\ln \alpha} > 0$, where $a > b > 0$.

**Proof.** See Appendix.

Proposition 3 shows that although it takes infinite time for $c°(t)$ and $k°(t)$ to approach the steady states, the needed time does differ: the larger is the difference between initial value and the stationary value, the longer time it takes to converge to the steady state. Hence, our model exhibits path-dependency, and postponement in the application of optimization practice has long-term effects. Specifically, the recovery time approaches infinity when $\alpha$ approaches 1, i.e., trivial difference in initial stock values may lead to irreversible outcome, at least in a time scale relevant to humanity. Our analysis thus shows that the potential costs of not taking the right action early may be very large; and the earlier effective action is taken, the less costly it will be.

**Proposition 4.** (Planning-horizon-variant capital, consumption, and shadow price) At each time $t$,

(i) the capital, $k_T(t)$, is monotonically increasing in the length of the planning horizon $T \in [0, \infty)$: $\frac{\partial k_T(t)}{\partial T} > 0$, for all $t \in [0, T]$;

(ii) the consumption, $c_T(t)$, is monotonically decreasing in $T$: $\frac{\partial c_T(t)}{\partial T} < 0$, for all $t \in [0, T]$;

(iii) the shadow price, $\lambda_T(t)$, is monotonically increasing in $T$: $\frac{\partial \lambda_T(t)}{\partial T} > 0$, for all $t \in [0, T]$.

**Proof.** See Appendix.

Proposition 4 shows that that capital, consumption, and shadow price all depend on the length of the planning horizon: the longer is the planning horizon, the higher are the capital and the accompanying shadow price, and the lower is the consumption at each time. Moreover, the shadow prices accompany an optimal path that extends to infinity should be higher than those with a finite planning horizon.

The above observations offer some intriguing insights, although the example we consider has been extremely simple. They seem to be in accordance with the traditional wisdoms on sustainability: that we should not discount the future, should improve technology, should make better plans concerning resource use, and should appreciate more the *in situ* value of the resources. What the sages of old perceived about sustainability by observing finite time has been extended to infinite time, precisely because we have extended their planning horizon to infinity without further modification.

## 4. Concluding Remarks

In this paper, we present a new approach to find optimums for problems in infinite-horizon dynamics with boundary conditions at infinite time without discounting. Interestingly, we show that most of the 'work' needed to obtain an optimum for undiscounted infinite planning horizon problems has already been done. One needs only to take the limit of the optimum obtained under the undiscounted finite planning horizon problem, if the model in consideration satisfies two fairly general conditions.

In particular, this approach is most relevant to the economics of climate change, in which the consensus on a zero time preference is forming. Since the ethical ramifications of discounting are nonexistent, the optimal paths for more complex models, obtained from our new approach, may exhibit some interesting new properties. We anticipate our assay to be a starting point for a systematic re-examination of infinite horizon optimization problems.

# APPENDIX

## PROOF OF LEMMA 1:

Let $\lim_{T \to \infty} b(T) = \beta$, then for an arbitrary $\varepsilon > 0$, there exists a $T_1 > 0$, such that if $T > T_1$, $\beta - \varepsilon < b(T) < \beta + \varepsilon$. We denote $\inf\{a(S) | S > T\}$ as $\underline{a}(T)$. Since $\lim_{T \to \infty} a(T) > 0$, there exists a $T_2$, such that if $T > T_2$, then $a(S) > 0$, that is, if $S > T_2$, then $a(S) > 0$.

Hence, if $S > \max(T_1, T_2)$, then $a(S)(\beta - \varepsilon) < (ab)(S) < a(S)(\beta + \varepsilon)$. Assuming that $T > \max(T_1, T_2)$, we have

$$\inf\{a(S)(\beta - \varepsilon) | S > T\} \leq \inf\{(ab)(S) | S > T\} \leq \inf\{a(S)(\beta + \varepsilon) | S > T\},$$

which implies

$$\underline{a}(T)(\beta - \varepsilon) \leq (\underline{ab})(T) \leq \underline{a}(T)(\beta + \varepsilon).$$

Hence,

$$\lim_{T \to \infty} \underline{a}(T)(\beta - \varepsilon) \leq \lim_{T \to \infty} (\underline{ab})(T) \leq \lim_{T \to \infty} \underline{a}(T)(\beta + \varepsilon).$$

Since $\varepsilon$ is arbitrarily chosen,

$$\lim_{T \to \infty} \underline{ab}(T) = \lim_{T \to \infty} \underline{a}(T) \beta. \qquad \text{Q.E.D.}$$

## PROOF OF PROPOSITION 3:

The difference of time needed for the two initial values, $k^a(0) \equiv \frac{1}{a}(\alpha)^{\frac{1}{1-\alpha}}$ and $k^b(0) \equiv \frac{1}{b}(\alpha)^{\frac{1}{1-\alpha}}$, to approach steady states is

$$\lim_{K^\circ(t) \to (1/\alpha)^{1/(1-\alpha)}} \left\{ \frac{1}{\ln \alpha} \ln\left( \frac{1}{\ln a}\left(\ln(\alpha)^{1/(1-\alpha)} - \ln K^\circ(t)\right)\right) - \frac{1}{\ln \alpha} \ln\left( \frac{1}{\ln b}\left(\ln(\alpha)^{1/(1-\alpha)} - \ln K^\circ(t)\right)\right) \right\}$$

$$= \frac{\ln\left(\frac{\ln a}{\ln b}\right)}{-\ln \alpha} > 0.  \qquad Q.E.D.$$

**PROOF OF PROPOSITION 4:**

(i) The proof proceeds by induction. We first show that $k_T(t)$ is monotonically increasing in $T$ when $t=1$. Since $k_T(1) = \left(1 - \frac{1-\alpha}{1-\alpha^{T+1}}\right) k(0)^\alpha$, it can be easily shown that $\frac{\partial k_T(1)}{\partial T} > 0$. Assume that $k_T(t)$ is monotonically increasing in $T$ at $t$, next we show that it is also valid at $t+1$. As $k_T(t+1) = \left(1 - \frac{1-\alpha}{1-\alpha^{T+1-t}}\right) k_T^\alpha(t)$, we have $\frac{\partial k_T(t+1)}{\partial T} > 0$. Hence, $k_T(t)$ is monotonically increasing in $T$ for all $t \in [0,T]$.

(ii) The proof proceeds also by induction. We first show that $c_T(t)$ is monotonically decreasing in $T$ when $t=1$. As $c_T(1) = \left(\frac{\alpha^{\alpha+1}(1-\alpha)}{(1-\alpha^T)^{1-\alpha}(1-\alpha^{T+1})^\alpha}\right) k(0)^{\alpha^2}$, it can be easily shown that $\frac{\partial c_T(1)}{\partial T} < 0$. Assume that $c_T(t)$ is monotonically decreasing in $T$ at $t$, next we show that it is also valid at $t+1$. Since $c_T(t) = \left(\frac{1-\alpha}{1-\alpha^{T+1-t}}\right) k_T^\alpha(t)$, we have

$$c_T(t+1) = \frac{1-\alpha}{1-\alpha^{T-t}}\left(\alpha\frac{1-\alpha^{T-t}}{1-\alpha}\cdot\frac{1-\alpha}{1-\alpha^{T+1-t}}k_T^{\alpha}(t)\right)^{\alpha} = (1-\alpha)^{1-\alpha}\alpha^{\alpha}\frac{(c_T(t))^{\alpha}}{(1-\alpha^{T-t})^{1-\alpha}}.$$

It can be easily shown that $\frac{\partial c_T(t+1)}{\partial T} < 0$. Hence, $c_T(t)$ is monotonically decreasing in $T$ for all $t \in [0,T]$.

(iii) Since $\lambda_T(t) = \frac{1}{c_T(t)}$, we immediately have $\frac{\partial \lambda_T(t)}{\partial T} > 0$ for all $t \in [0,T]$.      Q.E.D.

## PROOF OF THE FACT THAT THE PARAMETRIC EXAMPLE SATISFIES THE TWO CONDITIONS IN THEOREM 2:

(Here we present a formal proof that shows the two conditions are satisfied for the model we consider. For most examples, a numerical approach will be needed to verify the two conditions in Theorem 2. )

We first consider Condition (i) in Theorem 2. We have shown that $\ln c^{\circ}(t) < 0$, and

$$\sum_{t=0}^{\infty} \ln c^{\circ}(t) = -\infty.$$

$\forall M < 0$, $\exists N$ such that when $N \leq n$, we have

$$\sum_{t=0}^{n} \ln c^{\circ}(t) < M. \tag{A1}$$

We fix such an $n$. For all $0 \leq t \leq n$, since $\lim_{T\to\infty} c_T(t) = c^{\circ}(t)$, we see that

$$\lim_{T \to \infty} \left( \sum_{t=0}^{n} \ln c_T(t) \right) = \sum_{t=0}^{n} \ln c^{\circ}(t). \tag{A2}$$

From (A2), we see that $\exists T_0$, such that for $T_0 < T$, we have

$$\sum_{t=0}^{n} \ln c_T(t) < \sum_{t=0}^{n} \ln c^{\circ}(t) + 1. \tag{A3}$$

Let $T > n$, as $\ln c_T(t) < 0$, we have

$$\sum_{t=0}^{T} \ln c_T(t) = \sum_{t=0}^{n} \ln c_T(t) + \sum_{t=n+1}^{T} \ln c_T(t) < \sum_{t=0}^{n} \ln c_T(t). \tag{A4}$$

If $T > \max(n, T_0)$, from (A4), (A3), and (A1), we have

$$\sum_{t=0}^{T} \ln c_T(t) < \sum_{t=0}^{n} \ln c_T(t) < \sum_{t=0}^{T} \ln c^{\circ}(t) + 1 < M + 1,$$

And since $M$ is an arbitrarily chosen minus number, we see

$$\lim_{T \to \infty} \sum_{t=0}^{T} \ln c_T(t) = -\infty,$$

and the Condition (i) is satisfied.

We proceed to consider Condition (ii). We first show that $\ln(c_T(T))$ is bounded.

Note that

$$\ln(c_T(T)) = \ln(1-\alpha) + \frac{\alpha(1-\alpha^T)}{1-\alpha}\ln\alpha + \alpha^T \ln(k(0)) - (1-\alpha)\ln(1-\alpha)$$
$$-\frac{1-\alpha}{\alpha}\left[\alpha^2 \ln(1-\alpha^2) + \cdots + \alpha^T \ln(1-\alpha^T)\right] - \alpha^T \ln(1-\alpha^T)$$
$$= \frac{\alpha(1-\alpha^T)}{1-\alpha}\ln\alpha + \alpha^T \ln(k(0)) + \alpha \ln(1-\alpha)$$
$$-\frac{1-\alpha}{\alpha}\left[\alpha^2 \ln(1-\alpha^2) + \cdots + \alpha^T \ln(1-\alpha^T)\right] - \alpha^T \ln(1-\alpha^T).$$

As $-\alpha^t \ln(1-\alpha^t)$ is monotonically decreasing in $t$,

$$-\sum_{t=2}^{T} \alpha^t \ln(1-\alpha^t) \leq -\int_{1}^{T} \alpha^x \ln(1-\alpha^x) dx,$$

let $\alpha^x = y$,

$$-\int_{1}^{T} \alpha^x \ln(1-\alpha^x) dx = -\int_{\alpha}^{\alpha^T} \ln(1-y) \frac{dy}{\ln\alpha}$$

$$= \frac{1}{\ln\alpha}\left((1-\alpha^T)\ln(1-\alpha^T) - (1-\alpha^T) - (1-\alpha)\ln(1-\alpha) + (1-\alpha)\right),$$

we see the above is bounded by $\frac{1}{\ln\alpha}(-(1-\alpha)\ln(1-\alpha) - \alpha)$, implying that

$\lim_{T\to\infty}\left|\ln c_T(T)\right|$ is bounded by $\left|\frac{\alpha\ln\alpha}{1-\alpha} + \alpha\ln(1-\alpha)\right| + \frac{1-\alpha}{\alpha\ln\alpha}\left[-(1-\alpha)\ln(1-\alpha) - \alpha\right]$.

Let $G(\alpha) \equiv \frac{\alpha\ln\alpha}{1-\alpha} + \alpha\ln(1-\alpha) + \frac{1-\alpha}{\alpha\ln\alpha}\left[-(1-\alpha)\ln(1-\alpha) - \alpha\right]$, multiplying $G(\alpha)$ with $(1-\alpha)\alpha\ln\alpha < 0$, we see that

$$F(\alpha) \equiv G(\alpha)(1-\alpha)\alpha\ln\alpha$$
$$= \alpha^2(\ln\alpha)^2 + (1-\alpha)\alpha^2 \ln(1-\alpha)\ln\alpha - (1-\alpha)^3 \ln(1-\alpha) - (1-\alpha)^2 \alpha.$$

As $\ln(1+x) = x - \frac{1}{2}x^2 + \frac{1}{3}x^3 - \cdots$, $\ln(1-x) = -x - \frac{1}{2}x^2 - \frac{1}{3}x^3 - \cdots$, $F(\alpha)$ can be further restated as

$$F(\alpha) = \alpha^2 \left( -(1-\alpha) - \frac{1}{2}(1-\alpha)^2 - \frac{1}{3}(1-\alpha)^3 - \cdots \right)^2$$

$$+ (1-\alpha)\alpha^2 \left( -\alpha - \frac{1}{2}\alpha^2 - \frac{1}{3}\alpha^3 - \cdots \right)\left( -(1-\alpha) - \frac{1}{2}(1-\alpha)^2 - \cdots \right)$$

$$- (1-\alpha)^3 \left( -\alpha - \frac{1}{2}\alpha^2 - \frac{1}{3}\alpha^3 - \cdots \right) - (1-\alpha)^2 \alpha,$$

which is equivalent to

$$\frac{F(\alpha)}{\alpha} = \alpha \left( (1-\alpha) + \frac{1}{2}(1-\alpha)^2 + \frac{1}{3}(1-\alpha)^3 + \cdots \right)^2$$

$$+ (1-\alpha)\alpha \left( \alpha + \frac{1}{2}\alpha^2 + \frac{1}{3}\alpha^3 + \cdots \right)\left( (1-\alpha) + \frac{1}{2}(1-\alpha)^2 + \cdots \right)$$

$$+ (1-\alpha)^3 \left( 1 + \frac{1}{2}\alpha + \frac{1}{3}\alpha^2 + \cdots \right) - (1-\alpha)^2$$

$$= \alpha \left( (1-\alpha)^2 + (1-\alpha)^3 + \cdots \right)$$

$$+ (1-\alpha)\alpha \left( \alpha + \frac{1}{2}\alpha^2 + \frac{1}{3}\alpha^3 + \cdots \right)\left( (1-\alpha) + \frac{1}{2}(1-\alpha)^2 + \cdots \right)$$

$$+ (1-\alpha)^3 + (1-\alpha)^3 \left( \frac{1}{2}\alpha + \frac{1}{3}\alpha^2 + \cdots \right) - (1-\alpha)^2$$

$$= \alpha \left( (1-\alpha)^3 + \cdots \right)$$

$$+ (1-\alpha)\alpha \left( \alpha + \frac{1}{2}\alpha^2 + \frac{1}{3}\alpha^3 + \cdots \right)\left( (1-\alpha) + \frac{1}{2}(1-\alpha)^2 + \cdots \right)$$

$$+ (1-\alpha)^3 \left( \frac{1}{2}\alpha + \frac{1}{3}\alpha^2 + \cdots \right).$$

Hence, $F(\alpha) > 0$, $G(\alpha) < 0$, and $|\ln(c_T(T))|$ is bounded by

$$\left| \frac{\alpha \ln \alpha}{1-\alpha} + \alpha \ln(1-\alpha) \right| + \frac{1-\alpha}{\alpha \ln \alpha} \left[ -(1-\alpha)\ln(1-\alpha) - \alpha \right].$$

Next, we consider $\sum_{t=0}^{T-1} \ln(c_T(t)) - \sum_{t=0}^{T-1} \ln c^\circ(t)$. As

$$\sum_{t=0}^{T-1} \ln(c_T(t)) - \sum_{t=0}^{T-1} \ln c^\circ(t) = \sum_{t=0}^{T-1} \left( \ln(c_T(t)) - \ln c^\circ(t) \right),$$

where

$$\ln(c_T(t)) - \ln c^\circ(t) =$$
$$-(1-\alpha)\alpha^{-(T-t+1)} \left\{ \alpha^{T-t+1} \ln(1-\alpha^{T-t+1}) + \cdots + \alpha^T \ln(1-\alpha^T) \right\} - \alpha^t \ln(1-\alpha^{T+1}).$$

(A5)

Let $x = \alpha^s$, $T-t+1 \leq s \leq T$, the right-hand side of (A5) is

$$-(1-\alpha)\alpha^{-(T-t+1)} \sum_{s=T-t+1}^{T} \left( \alpha^s \ln(1-\alpha^s) \right) - \alpha^t \ln(1-\alpha^{T+1}),$$

the absolute value of which is bounded by the absolute value of

$$-(1-\alpha)\alpha^{-(T-t+1)} \int_{T-t}^{T} \alpha^s \ln(1-\alpha^s) ds - \alpha^t \ln(1-\alpha^{T+1}) \equiv (*).$$

Let $\alpha^s = x$, we then have

$$(*) = -(1-\alpha)\alpha^{-(T-t+1)} \int_{\alpha^{T-t}}^{\alpha^T} \ln(1-x) \frac{dx}{\ln \alpha} - \alpha^t \ln(1-\alpha^{T+1}).$$

Let $y = 1-x$, $(*)$ can be further rewritten as

$$\frac{(1-\alpha)\alpha^{-(T-t+1)}}{\ln \alpha} \int_{1-\alpha^{T-t}}^{1-\alpha^T} \ln y \, dy - \alpha^t \ln(1-\alpha^{T+1})$$

$$= \frac{(1-\alpha)\alpha^{-(T-t+1)}}{\ln\alpha}\left((1-\alpha^T)\ln(1-\alpha^T)+\alpha^T-(1-\alpha^{T-t})\ln(1-\alpha^{T-t})-\alpha^{T-t}\right)-\alpha^t\ln(1-\alpha^{T+1})$$

$$= \frac{(1-\alpha)}{\ln\alpha}\left((\alpha^{-(T-t+1)}-\alpha^{t-1})\ln(1-\alpha^T)+\alpha^{t-1}-(\alpha^{-(T-t+1)}-\alpha^{-1})\ln(1-\alpha^{T-t})-\alpha^{-1}-\frac{\ln\alpha}{1-\alpha}\alpha^t\ln(1-\alpha^{T+1})\right)$$

$$= -\frac{1-\alpha}{\ln\alpha}\left\{\underbrace{\left(\alpha^{-(T-t+1)}-\alpha^{-1}\right)\ln\left(1-\alpha^{(T-t)}\right)+\alpha^{-1}}_{(I)}-\underbrace{\left(\alpha^{-(T-t+1)}-\alpha^{t-1}\right)\ln\left(1-\alpha^T\right)}_{(II)}\underbrace{-\alpha^{t-1}}_{(III)}+\underbrace{\frac{\ln\alpha}{1-\alpha}\alpha^t\ln\left(1-\alpha^{T+1}\right)}_{(IV)}\right\}.$$

We observe that $\sum_{t=0}^{T-1}\left(\ln(c_T(t))-\ln c(t)\right)$ is bounded by the absolute value of

$$-\frac{1-\alpha}{\alpha\ln\alpha}\sum_{t=0}^{T-1}\left\{\underbrace{\left(\alpha^{-(T-t)}-1\right)\ln\left(1-\alpha^{(T-t)}\right)+1}_{(I)}-\underbrace{\left(\alpha^{-(T-t)}-\alpha^t\right)\ln\left(1-\alpha^T\right)}_{(II)}\underbrace{-\alpha^t}_{(III)}+\underbrace{\frac{\ln\alpha}{1-\alpha}\alpha^{t+1}\ln\left(1-\alpha^{T+1}\right)}_{(IV)}\right\}.$$

Next, we show that $\sum_{t=0}^{T-1}\left(\ln(c_T(t))-\ln c(t)\right)$ is bounded.

First, for term $(I)$, letting $s=T-t$, we can rewrite the sum of $(I)$ as

$$\sum_{s=1}^{T}\left\{(\alpha^{-s}-1)\ln(1-\alpha^s)+1\right\}.$$ Also, since $(\alpha^{-s}-1)\ln(1-\alpha^s)+1$ is positive, and

monotonically decreasing, we have

$$\sum_{s=1}^{T}\left\{(\alpha^{-s}-1)\ln(1-\alpha^s)+1\right\}$$
$$<(\alpha^{-1}-1)\ln(1-\alpha)+1+\int_{1}^{T}\left((\alpha^{-s}-1)\ln(1-\alpha^s)+1\right)ds.$$

Let $\alpha^s = x$, the above can be further rewritten as

$$(\alpha^{-1}-1)\ln(1-\alpha)+1+\int_1^T((\alpha^{-s}-1)\ln(1-\alpha^s)+1)ds$$

$$=(\alpha^{-1}-1)\ln(1-\alpha)+1+\int_\alpha^{\alpha^T}((x^{-1}-1)\ln(1-x)+1)\frac{dx}{(\ln\alpha)x}.$$

As $(\alpha^{-1}-1)\ln(1-\alpha)$ is bounded, and we focus on $\int_\alpha^{\alpha^T}((x^{-1}-1)\ln(1-x)+1)\frac{dx}{(\ln\alpha)x}$.

Let $h(x)=((x^{-1}-1)\ln(1-x)+1)x^{-1}$, we see that

$$\frac{dh(x)}{dx}=(-2x^{-3}+x^{-2})\ln(1-x)+(x^{-2}-x^{-1})\frac{-1}{1-x}$$

$$=\frac{-2+x}{x^3}\ln(1-x)-\frac{2}{x^2}$$

$$=\frac{-2+x}{x^3}\left(-x-\frac{1}{2}x^2-\frac{1}{3}x^3-\cdots\right)-\frac{2}{x^2}$$

$$=\frac{1}{x^3}\left(\frac{1}{6}x^3+\frac{2}{4\cdot 3}x^4+\frac{3}{5\cdot 4}x^5+\cdots+\frac{n-2}{n(n-1)}x^n+\cdots\right)>0.$$

Hence, we see that $h(x)$ is monotonically increasing in $x$, and we have

$$\lim_{x\to 0}h(x)=\lim_{x\to 0}((x^{-1}-1)\ln(1-x)+1)x^{-1}$$

$$=\lim_{x\to 0}\frac{-\ln(1-x)-x}{x^2}$$

$$=\lim_{x\to 0}\frac{\left(x+\frac{1}{2}x^2+\cdots\right)-x}{x^2}=\frac{1}{2}>0,$$

implying that $h(x)$ is bounded. Therefore, we have shown that sum of term ($I$) is bounded by $(\alpha^{-1}-1)\ln(1-\alpha)+1+\frac{1}{|\ln\alpha|}\int_0^\alpha((x^{-1}-1)\ln(1-x)+1)\frac{dx}{x}.$

For term ($II$), we see that

$$\sum_{t=0}^{T-1}\left(-\left(\alpha^{-(T-t)}-\alpha^{t}\right)\ln\left(1-\alpha^{T}\right)\right)$$

$$=-\left(\frac{\alpha^{-T}\left(1-\alpha^{T}\right)}{1-\alpha}-\frac{1-\alpha^{T}}{1-\alpha}\right)\ln\left(1-\alpha^{T}\right)$$

$$=\frac{-\alpha^{-T}+2-\alpha^{T}}{1-\alpha}\ln\left(1-\alpha^{T}\right)$$

$$=-\frac{\alpha^{-T}}{1-\alpha}\ln\left(1-\alpha^{T}\right)+\frac{2-\alpha^{T}}{1-\alpha}\ln\left(1-\alpha^{T}\right).$$

Obviously, the second term is small bounded. For the first term, as

$$\lim_{T\to\infty}\alpha^{-T}\ln\left(1-\alpha^{T}\right)=\lim_{x'\to 0}\frac{\ln\left(1-x'\right)}{x'}=\lim_{x'\to 0}\frac{\frac{-x'}{1-x'}}{x'}=-1,$$

we see that the sum of ($II$) is bounded by $\dfrac{\alpha}{1-\alpha}+1$.

For term ($III$), we see that $\sum_{t=0}^{T-1}\left(-\alpha^{t}\right)=-\dfrac{1-\alpha^{T}}{1-\alpha}$, hence, the sum of ($III$) is bounded by $-\dfrac{\alpha^{-1}}{1-\alpha}$.

Finally, for term ($IV$), we see that

$$\sum_{t=0}^{T-1}\frac{\ln\alpha}{1-\alpha}\alpha^{t+1}\ln\left(1-\alpha^{T+1}\right)=\frac{\ln\alpha}{1-\alpha}\alpha\frac{1-\alpha^{T}}{1-\alpha}\ln\left(1-\alpha^{T+1}\right),$$

which is small bounded. Hence, we have shown that $\sum_{t=0}^{T-1}\left(\ln\left(c_{T}(t)\right)-\ln c(t)\right)$ is bounded.

On the other hand, $\ln(k^{\circ}(T)^{\alpha}) = \ln\left(\alpha^{\frac{1-\alpha^{T}}{1-\alpha}} k(0)^{\alpha^{T}}\right)^{\alpha}$, the limit of which is $\ln(\alpha^{\frac{\alpha}{1-\alpha}})$.

Together with the fact that $\ln c_{T}(T)$ is bounded, we see that

$$-\left(\sum_{t=0}^{T-1}\ln\left(c_{T}(t)\right)-\sum_{t=0}^{T-1}\ln\left(c(t)\right)\right)-\ln c_{T}(T)+\ln\left(k(T)^{\alpha}\right)$$ is bounded, and as $\sum_{t=0}^{T}\ln\left(c_{T}(t)\right)$

is an infinite value, we see that Condition (ii) in Theorem 2 is also satisfied.

<div style="text-align: right">Q.E.D.</div>

**References**

Anand, S. and A. K. Sen (2000), "Human Development and Economic Sustainability", *World Development*, **28**, 2029-2049.

Anderson, R. M. (1991), "Chapter 39: Non-standard analysis with applications to economics", In Hildenbrand, W. & Sonnenschein, H. (Eds.) *Handbook of Mathematical Economics Vol. 4*, 2145-2208. North-Holland.

Arrow, K. J. (1999), "Discounting, Morality and Gaming", In Portney, P. R. & Weyant, J. P. (Eds.) *Discounting and Intergenerational Equity*. Washington, D.C.: Resources for the Future.

Brock, W.A. (1970), "On the Existence of Weakly Maximal Programmes in Multi-Sector Economy", *Review of Economic Studies*, **37**, 275-280.

Cline, W. R. (1992), *The Economics of Global Warming*. Washington, D.C.: Institute for International Economics.

Harrod, R. F. (1948), *Towards a Dynamic Economics*. London: Macmillan.

Heal, G. (1998), *Valuing the Future: Economic Theory and Sustainability*. New York: Columbia University.

Koopmans, T. (1960), "Stationary Ordinal Utility and Impatience", *Econometrica*, **28**, 287-309.

Nordhaus, W. D. (1994), *Managing the Global Commons: The Economics of Climate Change*. Cambridge: MIT Press.

Pigou, A. C. (1932), *The Economics of Welfare*. 4$^{th}$ ed. London: Macmillan.

Ramsey, F. P. (1928), "A Mathematical Theory of Saving", *Economic Journal*, **38**, 543-559.

Robinson, A. (1966), *Non-standard Analysis*. Amsterdam: North-Holland.

Rubio, J. E. (1994), *Optimization and Nonstandard Analysis*. Marcel Dekker, Inc.

Rubio, J. E. (2000), "Optimal Control Problems with Unbounded Constraint Sets", *Optimization*, **31**, 191-210.

Stern, N. (2007), *The Economics of Climate Change: The Stern Review*. Cambridge: Cambridge University Press.

Solow, R. M. (1974), "Richard T. Ely Lecture: The Economics of Resources or the Resources of Economics", *American Economic Review*, **64**, 1-14.

Solow, R. M. (1999), "Forward", In Portney, P. R. & Weyant, J. P. (Eds.) *Discounting and Intergenerational Equity*. Washington, D.C.: Resources for the Future.

Stokey, N. and Lucas R. (1989), *Recursive Methods in Economic Dynamics*. Cambridge: Harvard University Press.

Weitzman, M. L. (1999), "Just Keep Discounting, but…". In Portney, P. R. & Weyant, J. P. (Eds.) *Discounting and Intergenerational Equity.* Washington, D.C.: Resources for the Future.

Weitzman, M. L. (1998), "Why the Far-Distant Future Should Be Discounted at Its Lowest Possible Rate", *Journal of Environmental Economics and Management*, **36**, 201-208.